\documentclass[11pt]{amsart}

\usepackage{mathrsfs}
\usepackage{amsmath, amsfonts, amsthm, amssymb}
\usepackage[all,cmtip,line,arc]{xy}
\usepackage{array}
\usepackage{enumerate}
\usepackage{srcltx} 

\newtheorem{thm}{Theorem}[section]
\newtheorem{cor}[thm]{Corollary}
\newtheorem{prop}[thm]{Proposition}
\newtheorem{lem}[thm]{Lemma}

\theoremstyle{definition}
\newtheorem{defn}[thm]{Definition}

\theoremstyle{remark}
\newtheorem{rem}[thm]{Remark}

\newcommand \eps{\varepsilon}

\makeatletter
\let\c@equation\c@thm
\makeatother
\numberwithin{equation}{section}

\title{regularization of $m$-subharmonic functions and H\"older continuity}

\author{Jingrui Cheng, Yulun Xu}

\date{July 2022}

\begin{document}

\maketitle

\begin{abstract}
We use sup-convolution to find upper approximations of a bounded $m$-subharmonic function on a compact K\"ahler manifold with nonnegative holomorphic bisectional curvature.  As an application,  we show the H\"older continuity of solutions to $\sigma_m$ equation when the right hand side is in $L^p$,  $p>\frac{n}{m}$.  All these results generalize to more general complex Hessian equations.
\end{abstract}

\section{Introduction}
The regularization theory for plurisubharmonic functions on a compact complex manifold has been well-developed so far.  Let $M$ be a compact complex manifold and $\gamma$ a $(1,1)$ form on $M$.  Let $\varphi$ be a $\gamma$-plurisubharmonic function on $M$ (i.e.  $\gamma+\sqrt{-1}\partial\bar{\partial}\varphi\ge 0$ on $M$). 
Roughly speaking,  we hope to find a family of smooth functions $\varphi_{j}$ decreasing to $\varphi$ where each $\varphi_j$ is still $\gamma$-plurisubharmonic (or has very small loss of Hessian).  

Before we survey any previous result,  we would like to remark that classical version of such questions on $\mathbb{C}^n$ is quite easy.  Indeed,  if $\varphi$ is a plurisubharmonic function on $\mathbb{C}^n$,  we could simply define $\varphi_{\eps}(z)=\int_{\mathbb{C}^n}\varphi(z+\eps w)\chi(|w|)dV(w)$,  where $\chi$ is a standard smoothing kernel.  Then one can easily verify that $\varphi_{\eps}$ would decrease to $\varphi$ and $\varphi_{\eps}$ is also plurisubharmonic.
Such an approach even generalizes to the case when the underlying manifold $M$ is homogeneous.  Indeed,  it is possible to define ``translation with length $h$" on homogeneous manifolds and one can define $\varphi_{\eps}$ by taking the average over all translations of $\varphi$ with length $h$.  One may refer to \cite{EGZ3} for more details of this argument.

It is clear that we can no longer use the above definition for $\varphi_{\eps}$ on manifolds.  To find such approximations,  either we can try to define $\varphi_{\eps}$ locally using the above formula and try to glue them together,  or we can try to find a global expression to replace the above formula.  Blocki and Kolodziej \cite{BK} took the first approach for bounded plurisubharmonic functions.  They used convolution to define $\varphi_{\eps}$ locally,  and try to glue them together.  The major progress  was later made by Demailly \cite{Demailly},  \cite{DP}.  What he had in mind was mainly about applications to algebraic geometry.  Another application of his regularization theory was to prove the H\"older continuity of solutions to Monge-Ampere equations when the right hand side is in $L^p$ with $p>1$.

Let $(M,\omega_0)$ be a compact K\"ahler manifold.  Consider the complex Monge-Ampere equation:
\begin{equation*}
(\omega_0+\sqrt{-1}\partial\bar{\partial}\varphi)^n=e^F\omega_0^n,\,\,\sup_M\varphi=0.
\end{equation*}
We assume that $e^F\in L^p(\omega_0^n)$ for some $p>1$.  In the pioneering work by Kolodziej \cite{K1},  he proved that $\varphi$ has uniform $L^{\infty}$ bound (a bound for $\int e^F|F|^p\omega_0^n$ for $p>n$ is already sufficient).  
Later on,  Kolodziej \cite{K2} also proved that the solution to complex Monge-Ampere is H\"older continuous when the right hand side is in $L^p$ for some $p>1$.  Then S.  Dinew \cite{Dinew1} proved that the H\"older exponent actually depends only on $p$ and $n$,  under a positivity assumption on bisectional curvature.

In order to show that the solution is H\"older continuous,  the proof in \cite{DPGDHKZ} have two main ingredients (very roughly speaking):
\begin{enumerate}
\item (Stability Estimate)\\ Let $v$ be any bounded $\omega_0$-psh function,  then one has 
\begin{equation*}
\sup_M(v-\varphi)\le C||(v-\varphi)^+||_{L^1}^{\mu},
\end{equation*} 
for some $\mu>0$,  and $x^+=\max(x,0)$.
\item (Upper Regularization) \\Consider $\varphi_{\eps}(z)=\int_{\zeta\in T_zM}\varphi(\exp_z(\eps \zeta))\chi(|\zeta|^2)dV_z(\zeta)$,  and show that $\omega_0+\sqrt{-1}\partial\bar{\partial}\varphi_{\eps}\ge -\eps^{\alpha}\omega_0$,  $||\varphi_{\eps}-\varphi||_{L^1}\le C\eps^{\beta}$.  
\end{enumerate}
Once we have (1) and (2),  we can take $v=\frac{\varphi_{\eps}}{1+\eps^{\alpha}}$ in the Stability Estimate,  and obtain that $\varphi_{\eps}-\varphi\le C\eps^{\gamma}$.  On the other hand,  $\Delta \varphi\ge -n$ since $\varphi$ is $\omega_0$-psh,  we will be able to get $\sup_{z'\in B_{\frac{\eps}{2}}(z)}\varphi(z')-\varphi(z)\le C\eps^{\gamma}$,  and this will imply $\varphi$ is in $C^{\gamma}$.

The goal of this present work is to make a first attempt to generalize such results to more general Hessian equations.  We will consider equations of the following form:
\begin{equation}\label{1.1New}
f(\lambda[h_{\varphi}])=e^{F},\,\,\lambda[h_{\varphi}]\in\Gamma.
\end{equation}
In the above,  $\varphi$ is a real valued function on $M$ and $(h_{\varphi})_j^i=\sum_kg^{i\bar{k}}(g_{j\bar{k}}+\varphi_{j\bar{k}})$ and $\lambda[h_{\varphi}]$ means the eigenvalues of $\{(h_{\varphi})_j^i\}_{1\le i,\,j\le n}$.
In the above,  $f(\lambda_1,\lambda_2,\cdots,\lambda_n)$ is a function defined on a cone $\Gamma\subset \mathbb{R}^n$.  Following \cite{GPT},  we make the following assumptions on the cone $\Gamma$ as well as the function $f$:
\begin{enumerate}
\item $\Gamma$ and $f$ are symmetric in $\lambda$,  which means they are invariant under swapping $\lambda_i$ and $\lambda_j$.
\item $\Gamma_n\subset\Gamma\subset\Gamma_+$,  where $\Gamma_n=\{\lambda\in \mathbb{R}^n:\lambda_i>0,\,1\le i\le n\}$ and $\Gamma_+=\{\lambda\in\mathbb{R}^n:\sum_i\lambda_i>0\}$.  $\Gamma$ is convex.
\item $\frac{\partial f}{\partial \lambda_i}>0$,  $1\le i\le n$ and $f$ is concave on $\Gamma$.
\item There exists $c_0>0$ such that if we put $F(h)=f(\lambda(h))$ where $h$ is a positive definite $n\times n$ Hermitian matrix,  then $\det \frac{\partial F}{\partial h_{i\bar{j}}}\ge c_0>0$.
\item There exists $C_0>0$ such that for $\lambda\in\Gamma_n$,  $\sum_{i=1}^n\lambda_i\frac{\partial f}{\partial\lambda_i}\le C_0f$.
\end{enumerate}
The examples which one would like to study first is the $\sigma_m$ equation (with $1\le m\le n$),  where $f(\lambda)=\sigma_m^{\frac{1}{m}}(\lambda)$ where \sloppy $\sigma_m(\lambda)=\sum_{1\le i_1<i_2<\cdots <i_m\le n}\lambda_{i_1}\cdots \lambda_{i_m}$,  defined on the cone $\Gamma_m=\{\lambda\in\mathbb{R}^n:\sigma_i(\lambda)>0,\,1\le i\le m\}$.
Note that examples like $f(\lambda)=\big(\frac{\sigma_k}{\sigma_m}\big)^{\frac{1}{k-m}}$,  $1\le m\le k-1$ would not satisfy the assumption (4) above,  but $f(\lambda)=\big(\frac{\sigma_k}{\sigma_m}\big)^{\frac{1}{k-m}}+c\sigma_l^{\frac{1}{l}}$ with $c>0$ does satisfy the assumption (4).

Due to the main result of \cite{GPT},  $\varphi$ is bounded in $L^{\infty}$.  In order to study the H\"older continuity,  we need to establish the Stability Result and Upper Regularization result for (\ref{1.1New}),  similar to the approach for complex Monge-Ampere equations.

The stability result for a general Hessian equations can be obtained following a PDE approach developed in \cite{WWZ0},  \cite{WWZ} and \cite{GPT}.  It is done in Section \ref{Stability}.  It can be stated as follows:
\begin{prop}\label{p1.2}
(Lemma \ref{l3.7}) Let $\delta>0$ and let $v$ be a bounded function such that $\lambda\big(\frac{\delta}{2}I+(h_v)_j^i\big)\in \Gamma$.  Let $\varphi$ solve \ref{1.1New} with $e^{nF}\in L^{p_0}(\omega_0^n)$ for some $p_0>1$.  Let $s_0>0$ be suitably large (compared with $\delta$) such that
\begin{enumerate}
\item $s_0\ge 2\delta||v||_{L^{\infty}}$;
\item $\int_M((1-\delta)v-\varphi-s_0)^+e^{nF}\omega_0^n\le \delta^{n+1}$.
\end{enumerate}
Then for any $\mu<\frac{1}{nq_0}$ (where $q_0=\frac{p_0}{p_0-1}$),  
\begin{equation*}
\sup_M(v-\varphi)\le s_0+C_1s_0^{-\mu}||(v-\varphi)^+||_{L^1}^{\mu}.
\end{equation*}
\end{prop}

The next step is to find a suitable upper approximation $\varphi_{\eps}$.  We need to choose $\varphi_{\eps}$ so that: 
\begin{enumerate}
\item ($\varphi_{\eps}$ is almost $\Gamma$-admissible)\\There exists $\gamma_1>0$ such that $\frac{C\eps^{\gamma_1}}{2}I+(h_{\varphi_{\eps}})_j^i\in \Gamma$ in the viscosity sense.  (Actually one can show that $(\varphi_{\eps})_{i\bar{j}}$ is bounded below for each $\eps>0$)
\item ($\varphi_{\eps}$ approximates $\varphi$ in $L^1$) \\There exists $\gamma_2>0$ such that $||\varphi_{\eps}-\varphi||_{L^1}\le C\eps^{\gamma_2}$.
\item ($\varphi_{\eps}$ approximates $\varphi$ from above) \\$\varphi_{\eps}(z)\ge \sup_{z'\in B_{\eps}(z)}\varphi(z')-C\eps$.
\end{enumerate}
In this work,  we try a sup convolution to construct the approximation,  namely we define
\begin{equation}\label{1.3N}
\varphi_{\eps}(z)=\sup_{\xi\in T_zM}\big(\varphi(\exp_z(\xi))+\eps-\frac{1}{\eps}|\xi|_z^2)\big).
\end{equation}
In the above,  $\exp_z(\xi)$ is the exponential map defined using $\omega_0$,  and $|\xi|_z$ is the length of the tangent vector (at $z$).  
Unfortunately,  at this moment,  we need to require $M$ to have nonnegative holomorphic bisectional curvature.  With this assumption,  one has:
\begin{prop}\label{p1.4}
(Theorem \ref{hessian L1 estimate}) Assume that $M$ has nonnegative holomorphic bisectional curvature.  Let $\varphi$ be $\Gamma$-admissible and bounded in the viscosity sense.  Define $\varphi_{\eps}$ according to (\ref{1.3N}),  then (1),  (2),  (3) above holds.
\end{prop}
The requirement of being compact K\"ahler and having nonnegative holomorphic bisectional curvature is rather restrictive.  Mok \cite{Mok} has a uniformization theorem which says that its universal cover is isometrically biholomorphic to $(\mathbb{C}^k,g_0)\times (\mathbb{P}^{N_1},\theta_1)\times\cdots(\mathbb{P}^{N_l},\theta_l)\times (M_1,g_1)\times\cdots(M_p,g_p)$.
In the above,  $g_0$ is the standard Euclidean metric.  $M_1$,  $\cdots$,  $M_p$ are Hermitian symmetric spaces with $g_1,\cdots g_p$ being canonical metrics on them,  and $\theta_1$,  $\cdots$,  $\theta_l$ are metrics on $\mathbb{P}^{N_1}$,  $\cdots$,  $\mathbb{P}^{N_l}$ carry nonnegative holomorphic bisectional curvature.
In our work,  the only reason we need this additional curvature assumption is to estimate the lower bound of complex Hessian for $\varphi_{\eps}$ (i.e.  to show $\varphi_{\eps}$ is almost $\Gamma$-admissible in the viscosity sense).
We hope to relax or remove such curvature assumptions in the future.

Combining Proposition \ref{p1.2} and \ref{p1.4},  we can show that solutions to (\ref{1.1New}) with $e^F\in L^{p_0}$ for some $p_0>1$ is H\"older continuous.  Indeed,  we take $v=\frac{1}{1+\frac{C\eps^{\gamma_1}}{2}}\varphi_{\eps}$,  then $v$ would satisfy the assumptions in Proposition \ref{p1.2}.  Then we can take $\delta=C\eps^{\gamma_1}$ and we can take $s_0=\eps^a$ for some $a>0$ and small.  This is possible due to some additional estimates regarding (1) and (2) in Proposition \ref{p1.2}.  So the estimate would give us that:
\begin{equation*}
\sup_M(v-\varphi)\le C\eps^b.
\end{equation*}
Since $\varphi_{\eps}$ is uniformly bounded,  this easily translates to:
\begin{equation*}
\sup_{B_{\eps}(z)}\varphi(z')\le \sup_M(\varphi_{\eps}-\varphi)+C\eps\le C'\eps^{b'}.
\end{equation*}
This would imply that $\varphi$ is H\"older continuous.
More precisely,  we have:
\begin{thm}
(Theorem \ref{t4.1})
Let $(M,\omega_0)$ be a compact K\"ahler manifold with nonnegative holomorphic bisectional curvature.  Let $\varphi$ be a solution to (\ref{1.1New}) with $f$ satisfying the structural assumptions and the right hand side $e^F\in L^{p}(\omega_0^n)$ for some $p>n$.  Then for any $\gamma<\frac{p-n}{2(p-n+pn)}$ we have that $||\varphi||_{C^{\gamma}}\le C'$,  with $C'$ depends on the background metric,  $p_0$,  $n$,  as well as $||e^F||_{L^{p}}$ as well as the choice of $\gamma$.
\end{thm}
This theorem corresponds to Theorem \ref{t4.1},  with $p_0=\frac{p}{n}$.

If we specialize to the $\sigma_m$ equation by taking $f(\lambda)=\sigma_m^{\frac{1}{m}}(\lambda)$,  then we have:
\begin{cor}
Let $(M,\omega_0)$ be a compact K\"ahler manifold with nonnegative holomorphic bisectional curvature.  Let $\varphi$ be an $m$-subharmonic solution to the $\sigma_m$ equation:
\begin{equation*}
(\omega_0+\sqrt{-1}\partial\bar{\partial}\varphi)^m\wedge \omega_0^{n-m}=e^F\omega_0^n,\,\,\,\lambda\big(g^{i\bar{k}}(g_{\varphi})_{j\bar{k}}\big)_j^i\in\Gamma_m.
\end{equation*}
Assume that $e^F\in L^p(\omega_0^n)$ for some $p>\frac{n}{m}$,  then for any $\gamma<\frac{mp-n}{2(mp-n+mpn)}$,  one has $||\varphi||_{C^{\gamma}}\le C'$,  where $C'$ depends additionally on the background metric and $||e^F||_{L^p(\omega_0^n)}$ and also the choice of $\gamma$.
\end{cor}
As far as we know,  the previous results regarding H\"older continuity to solutions of $\sigma_m$-equation are for domains of $\mathbb{C}^n$ and $\omega_0$ is the standard Euclidean metric in $\mathbb{C}^n$,  see \cite{BAZ} and \cite{NC}.
The organization of the paper is as follows.

In Section 2,  we develop the theory for upper approximation for bounded functions for which the comlex Hessian lies in a general cone $\Gamma$,  assuming the holomorphic bisectional curvature of the underlying manifold is nonnegative.

In Section 3,  we use the approximations found in section 2 to show the H\"older continuity of solutions to complex Hessian equations.

\section{Approximation using sup-convolution}
In this section,  we verify that the $\varphi_{\eps}$ defined by (\ref{1.3N}) will be a good upper approximation of $\varphi$,  satisfying (1),  (2) and (3) above (\ref{1.3N}).  More precisely,  we have the following theorem:

\begin{thm}\label{hessian L1 estimate}
Suppose that $(M,\omega_0)$ is a K\"ahler manifold such that the holomorphic bisectional curvature is non-negative.  Let $\varphi\in C(M)$ and $\lambda[h_{\varphi}]\in \Gamma$ in the viscosity sense,  as defined by Definition \ref{d2.3}.
Define $\varphi_{\eps}$ as in (\ref{1.3N}).  Then we have:\\
$(i)$ $\varphi_{\eps}$ is semi-convex and $|\varphi_{\eps}|\le C$ for a constant $C$ depending only on $||\varphi||_{L^{\infty}}$. \\
$(ii)$ There exists a constant $C$ depending on $||\varphi||_{\infty}$ and $(M,\omega_0)$ such that:
\begin{equation*}
    \lambda(g^{i\bar{k}} (\varphi_{\eps})_{i\bar{k}} +(1+C\eps^{\frac{1}{2}})I)\in \Gamma.
\end{equation*}
$(iii)$ There exists a constant $C$ depending on $||\varphi||_{\infty}$ and $(M,\omega_0)$ such that:
\begin{equation*}
    ||\varphi_{\eps}-\varphi||_{L^1}\le C \eps^{\frac{1}{4}}.
\end{equation*}
Moreover, if we have that $\varphi$ is H\"older continuous with a H\"older exponent $\gamma_0$, then we can improve the above results as follows:\\
$(ii')$ There exists a constant $C$ depending on $||\varphi||_{C^{0,\gamma_0}}$ and $(M,\omega_0)$ such that:
\begin{equation*}
    \lambda(dd^c \varphi_{\eps} +(1+C\eps^{\frac{1+\gamma_0}{2-\gamma_0}})\omega_0)\in \Gamma.
\end{equation*}
$(iii')$ There exists a constant $C$ depending on $||\varphi||_{C^{0,\gamma_0}}$ and $(M,\omega_0)$ such that:
\begin{equation*}
    ||\varphi_{\eps}-\varphi||_{L^1}\le C \eps^{\frac{1}{2(2-\gamma_0)}}.
\end{equation*}
\end{thm}

Demaily used a special holomorphic coordinate in his paper \cite{DP} to simplify the calculation. We give the following definition in light of his idea.
\begin{defn}\label{d2.2}
We call a holomorphic coordinate $(z)$ normalized coordinate if in this coordinate the K\"ahler form has the following expansion:
\begin{equation*}
    (\omega_0)_{i\bar{k}}=\delta_{i k} -c_{i k \alpha \beta} z_{\alpha} \bar{z}_{\beta}+O(|z|^3)
\end{equation*}
\end{defn}

In \cite{DP},  Demailly define the following approximation for an $\omega_0$-psh function $\varphi$:
\begin{equation*}
\Psi(z,w)=\int_{\zeta\in T_zM}\psi(exph_z(w\zeta))\chi(|\zeta|^2)d\lambda(\zeta).
\end{equation*}
In the above,  $\chi$ is a smoothing kernel with compact support.  exph is a modification of the usual exponential map which erases all the non-holomorphic terms with respect to $\zeta$ in the Taylor expansion of $\exp_z(\zeta)$ at $\zeta=0$.  The estimates in \cite{DP} depends crucially on that $\sqrt{-1}\partial\bar{\partial}\psi\ge 0$,  and does not seem to generalize to when it may have negative eigenvalues.

Now we work with (\ref{1.3N}) and our first step is to get a lower bound of the complex Hessian.  In this work,  we will mostly work with the notion of viscosity solutions.  Following \cite{JLS},  we can define:
\begin{defn}\label{d2.3}
Let $\varphi\in C(M)$,  we say that $\lambda[h_{\varphi}]\in\Gamma$ in the viscosity sense,  if for all $x_0\in M$ and all $P\in C^2(M)$ touching $\varphi$ from above at $x_0$ ($P(x_0)=\varphi(x_0)$,  $P(x)\ge \varphi(x)$ in a neighborhood of $x_0$),  we have
\begin{equation*}
\lambda[h_P](x_0)=\lambda\big((g^{i\bar{k}}(g_{j\bar{k}}+\partial_j\partial_{\bar{k}}P))_j^i\big)(x_0)\in \Gamma.
\end{equation*}
\end{defn}
\begin{rem}
We note that if $\varphi\in C^2(M)$ with $\lambda[h_{\varphi}]\in\Gamma$ pointwise,  then $\lambda[h_{\varphi}]\in\Gamma$ also in viscosity sense.  Indeed,  if $P$ touches $\varphi$ from above at $x_0$,  then we would have $P_{i\bar{j}}(x_0)\ge \varphi_{i\bar{j}}(x_0)$.  Hence if we take eigenvalues,  $\lambda_i(\lambda[h_P])\ge \lambda_i(\lambda[h_{\varphi}])$,  $1\le i\le n$.  But $\Gamma$ is a convex cone containing $\Gamma_n:=\{\lambda\in\mathbb{R}^n:\lambda_i>0\}$,  we get $\lambda[h_P]\in\Gamma$.
\end{rem}

First we show that $\lambda[h_{\varphi}]\in\Gamma$ in the viscosity sense translate to a similar condition for $\varphi_{\eps}$.  We start with the following simple observation.
\begin{lem}\label{l2.5}
\begin{enumerate}
\item
Let $\varphi\in C(M)$ and we define $\varphi_{\eps}$ according to (\ref{1.3N}).  Let $x_0\in M$ and $P$ be a $C^2$ function defined in a neighborhood of $x_0$ touching $\varphi_{\eps}$ from above.  Assume that $\varphi_{\eps}(x_0)=\varphi(\exp_{x_0}(\xi_0))+\eps-\frac{1}{\eps}|\xi_0|^2_{x_0}$ for some $\xi_0\in T_{x_0}(M)$.  Let $\xi(z)\in T_zM$ be a smooth vector field defined in a neighborhood of $x_0$ with $\xi(x_0)=\xi_0$.  Define $\phi(z)=\exp_z(\xi(z))$,  then $z\mapsto P-\eps+\frac{1}{\eps}|\xi(z)|^2_z$ touches $\varphi\circ \phi$ from above at $z_0$.
\item If the map $\phi$ is holomorphic with $D_z\phi(x_0)$ being invertible,  then $w\mapsto (P-\eps+\frac{1}{\eps}|\xi|^2)\circ \phi^{-1}(w)$ touches $\varphi$ from above at $\exp_{x_0}(\xi_0)$.
\end{enumerate}
\end{lem}
\begin{proof}
The proof is rather straightforward out of the definition.  For part (1),  let $P$ touch $\varphi_{\eps}$ from above will lead to the following string of inequalities:
\begin{equation*}
P(z)\ge \varphi_{\eps}(z)\ge \varphi(\exp_z(\xi(z)))+\eps-\frac{1}{\eps}|\xi(z)|^2_z,
\end{equation*}
with equality achieved at $x_0$.
This just means $P(z)+\frac{1}{\eps}|\xi(z)|^2_z-\eps$ touches $\varphi\circ \phi$ from above.

For part (2),  the additional assumption on $\phi$ implies that $\phi$ is actually biholomorphic between a neighborhood of $x_0$ and a neighborhood of $w_0:=\exp_{x_0}(\xi_0)$.
\end{proof}

As a direct consequence,  we have
\begin{cor}\label{c2.6}
Assume that $\lambda[h_{\varphi}]\in\Gamma$ in the viscosity sense.  Let $\xi(z)$ and $\phi(z)$ satisfy the assumptions in point (2) of Lemma \ref{l2.5}.  Put $w_0=\exp_{z_0}(\xi_0)$ and $Q(w)=(P-\eps+\frac{1}{\eps}|\xi(z)|^2_z)\circ \phi^{-1}(w)$.  Then $\lambda(g^{i\bar{k}}(g_{j\bar{k}}+\partial_{w_j\bar{w}_k}Q)(w_0)\in \Gamma$.
\end{cor}

In order to get the estimate of Hessian for $\varphi_{\eps}$ (in the viscosity sense),  we need to translate Corollary \ref{c2.6} to a statement at $z_0$.  For this we have:
\begin{lem}\label{l2.7}
$\lambda(g^{i\bar{k}}(g_{j\bar{k}}+\partial_{w_j\bar{w}_k}Q))(w_0)$ is the same as the set of roots of
\begin{equation*}
f(\lambda)=\det \big((\lambda-1) \sum_{i,j}g_{i\bar{j}}(w_0)\frac{\partial \phi_i}{\partial z_a}\frac{\partial \phi_{\bar{j}}}{\partial\bar{z}_b}(x_0)-(P+\eps-\frac{1}{\eps}|\xi(z)|^2_z)_{z_a\bar{z}_b}(x_0)\big)_{1\le a,b\le n}.
\end{equation*}
\end{lem}
\begin{proof}
First we observe that
\begin{equation}\label{2.8NN}
Q_{i\bar{j}}(w_0)=\frac{\partial(\phi^{-1})_a}{\partial w_i}\frac{\partial(\phi^{-1})_{\bar{b}}}{\partial\bar{w}_j}(P-\eps+\frac{1}{\eps}|\xi(z)|^2_z)_{z_a\bar{z}_b}(x_0).
\end{equation}
On the other hand,  $\lambda(g^{i\bar{k}}(g_{j\bar{k}}+\partial_{w_j\bar{w}_k}Q))(w_0)$ is the set of roots for:
\begin{equation*}
\begin{split}
&\det\big(\lambda \delta_{ij}-g^{i\bar{k}}(g_{j\bar{k}}+\partial_{j\bar{k}}Q)(w_0)\big)=\det\big((\lambda-1) \delta_{ij}-g^{i\bar{k}}\partial_{j\bar{k}}Q(w_0)\big)\\
&=\det(g^{-1})\det\big((\lambda-1)g_{i\bar{j}}(w_0)-Q_{i\bar{j}}(w_0)\big)\\
&=\det(g^{-1})|\det(D_w\phi^{-1})(w_0)|^2\det\big((\lambda-1)\sum_{i,j}g_{i\bar{j}}(w_0)\frac{\partial\phi_i}{\partial z_a}\frac{\partial\phi_{\bar{j}}}{\partial\bar{z}_b}(z_0)\\
&-(P+\eps-\frac{1}{\eps}|\xi(z)|^2_z)_{z_a\bar{z}_b}\big).
\end{split}
\end{equation*}
In the last equality,  we used (\ref{2.8NN}).
\end{proof}
Since we want to ``translate" from $w_0=\exp_{z_0}(\xi_0)$,  we wish to take the map $\phi$ to satisfy:
\begin{equation}\label{2.9N}
\sum_{i,j}g_{i\bar{j}}(w_0)\frac{\partial\phi_i}{\partial z_a}\frac{\partial\phi_{\bar{j}}}{\partial\bar{z}_b}(z_0)=g_{a\bar{b}}(x_0),\,\,\phi \text{ holomorphic},\,\,\det D_z\phi(x_0)\neq 0. 
\end{equation}
Next let us explain how to achieve (\ref{2.9N}).  First let us observe that $w_0$ and $x_0$ are actually very close if $\eps$ is small enough.  Indeed,
\begin{lem}\label{l2.10}
Let $\varphi\in C(M)$ and $\varphi_{\eps}$ be as defined by (\ref{1.3N}).  Let $\xi_0\in M$ and $\xi_0\in T_{x_0}M$.  We assume that $\xi_0$ achieves the sup in the definition of $\varphi_{\eps}$,  then 
\begin{equation*}
|\xi_0|_{x_0}\le (2||\varphi||_{L^{\infty}}\eps)^{\frac{1}{2}}.
\end{equation*}
In particular,  $d(w_0,x_0)\le C\eps^{\frac{1}{2}}$.  

If we assume additionally that $\varphi\in C^{\alpha}(M)$,  then 
\begin{equation*}
|\xi_0|_{x_0}\le (\eps[\varphi]_{\alpha})^{\frac{1}{2-\alpha}},
\end{equation*}
where $[\varphi]_{\alpha}=\sup_{x,y\in M}\frac{|\varphi(x)-\varphi(y)|}{d(x,y)^{\alpha}}$.  
\end{lem}
\begin{proof}
Since $\xi_0$ achieves the sup,  we have:
\begin{equation*}
\varphi(\exp_{x_0}(\xi_0))+\eps-\frac{1}{\eps}|\xi_0|^2_{x_0}\ge \varphi(x_0)+\eps-0.
\end{equation*}
In the right hand side above,  we are taking $\xi=0$.  If $\varphi$ does not have any estimate on the modulus of continuity,  we can only estimate $\varphi(\exp_{x_0}(\xi_0))-\varphi(x_0)$ by $2||\varphi||_{L^{\infty}}$.  If $\varphi\in C^{\alpha}$,  then we have
\begin{equation*}
\varphi(\exp_{x_0}(\xi_0))-\varphi(x_0)\le [\varphi]_{\alpha}d(\exp_{x_0}(\xi_0),x_0)^{\alpha}\le [\varphi]_{\alpha}|\xi_0|_{x_0}^{\alpha}.
\end{equation*}
\end{proof}
As a direct consequence of this,  we note that $\varphi_{\eps}$ is semi-convex:
\begin{lem}\label{l2.11New}
Define $\varphi_{\eps}$ according to (\ref{1.3N}).  Then $\varphi_{\eps}$ is semi-convex for $\eps$ small enough,  in the sense that for any $z_0\in M$,  there is a coordinate neighborhood of $z_0$ such that $\varphi_{\eps}(z)+C_{\eps}|z-z_0|^2$ is locally convex near $z_0$.
In particular,  $\varphi_{\eps}$ is second order differentiable a.e.
\end{lem}
\begin{proof}
Since for any $z\in M$,  $D_{\xi}\exp_z(\xi)|_{\xi=0}$ is invertible,  we see that,  by implicity function theorem,  there exists a neighborhood $U_0$ of $z_0$,  such that for any $z\in U_0,\,w\in U_0$,  there is a unique $\xi\in T_zM$ such that $\exp_z(\xi)=w$.  We can denote $\xi=\log_z(w)$ and we may assume that this map is smooth.

Because of Lemma \ref{l2.10},  we have
\begin{equation*}
\varphi_{\eps}(z)=\sup_{|\xi|_z\le (2||\varphi||_{L^{\infty}})^{\frac{1}{2}}}(\varphi(\exp_z\xi)+\eps-\frac{1}{\eps}|\xi|^2)=\sup_{w\in U_0}\big(\varphi(w)+\eps-\frac{1}{\eps}|\log_z(w)|^2_z\big).
\end{equation*}
We observe that there exists $C_{\eps}>0$ such that 
\begin{equation*}
D_z^2\big(\varphi(w)+\eps-\frac{1}{\eps}|\log_z(w)|^2_z\big)\ge -C_{\eps}I,\,\,\,\text{for any $w\in U_0$.}
\end{equation*}
This essentially follows from that $(z,w)\mapsto -\frac{1}{\eps}|\log_zw|^2_z$ is jointly smooth.
Therefore,  taking the sup with respect to $w$ will make it semi-convex in $z$ with the same lower bound of the Hessian.
\end{proof}
Next we are going to take the normal coordinate centered at $x_0$ as in Definition \ref{d2.2}.  Without loss of generality,  we can assume that $x_0$ is represented by $z=0$ under this coordinate.  Let $N$ be an $n\times n$ matrix such that:
\begin{equation}\label{2.11N}
\sum_{i,j}N_{ia}\bar{N}_{jb}g_{i\bar{j}}(w_0)=\delta_{ab}=g_{a\bar{b}}(x_0).
\end{equation}
Then under this coordinate,  we define 
\begin{equation}\label{2.11}
\phi(z)=\exp_{x_0}(\xi_0)+N\cdot z.
\end{equation}
We need to make sure that one can define $\xi(z)$ which satisfies $\exp_z(\xi(z))=\phi(z)$:
\begin{lem}\label{l2.13}
For all $\eps>0$ small enough,  there exists a smooth vector field $z\mapsto \xi(z)\in T_zM$ defined in a neighborhood of $x_0$ with $\xi(x_0)=\xi_0$ such that $\exp_z(\xi(z))=\phi(z)$.
\end{lem}
\begin{proof}
This follows from the implicit function theorem.
Indeed,  let $U_0$ and $V_0$ be neighborhoods of $0$ and $\xi_0$ respectively (in $\mathbb{C}^n$),  then we may consider the nonlinear map (where the exponential map is represented using the above normal coordinates):
\begin{equation*}
F:U_0\times V_0\rightarrow \mathbb{C}^n,\,\,\,(z,\xi)\mapsto \exp_z(\xi)-\exp_{x_0}(\xi_0)-N\cdot z.
\end{equation*}
It is clear that $F(0,\xi_0)=0$.  Also we find that $D_{\xi}F(0,\xi_0)=D_{\xi}(\exp_{x_0}(\xi))|_{\xi=\xi_0}$.  Note that $D_{\xi}(\exp_{x_0}(\xi))|_{\xi=0}=id$,  we see that $D_{\xi}F(0,\xi_0)$ will be invertible if $\eps$ is small enough,  thanks to Lemma \ref{l2.10}.  Then we can apply the implicity function theorem to define $\xi(z)$ with $F(z,\xi(z))=0$.
\end{proof}

Combining Corollary \ref{c2.6} and Lemma \ref{l2.7},  we can conclude that:
\begin{prop}\label{p2.14}
Let $\varphi\in C(M)$ and $\lambda[h_{\varphi}]\in \Gamma$ in the viscosity sense.  Define $\varphi_{\eps}$ according to (\ref{1.3N}).  Let $\phi$ be given by (\ref{2.11}) and $\xi(z)$ be given by Lemma \ref{l2.13}.  Let $P\in C^2(M)$ and touch $\varphi_{\eps}$ from above at $x_0\in M$,  then 
\begin{equation*}
\lambda\big(g^{i\bar{k}}(g_{j\bar{k}}+\partial_{z_j\bar{z}_k}P+(\frac{1}{\eps}|\xi(z)|^2_z)_{z_j\bar{z}_k})(x_0)\big)\in\Gamma.
\end{equation*}
\end{prop}
At this moment,  it is clear that we have to estimate the lower bound of the complex Hessian for $\frac{1}{\eps}|\xi(z)|^2_z$ at $x_0$.
We will do it now.  
First,
\begin{equation}\label{2.15N}
\begin{split}
&(|\xi(z)|^2)_{i\bar{j}}(x_0)=(g_{a\bar{b}}(z)\xi_a(z)\bar{\xi}_{\bar{b}}(z))_{i\bar{j}}=g_{a\bar{b},i\bar{j}}(x_0)\xi_a\bar{\xi}_{\bar{b}}(x_0)+g_{a\bar{b}}(x_0)(\xi_a\bar{\xi}_{\bar{b}})_{i\bar{j}}(x_0)\\
&+g_{a\bar{b},i}(x_0)(\xi_a\bar{\xi}_{\bar{b}})_{\bar{j}}(x_0)+g_{a\bar{b},\bar{j}}(x_0)(\xi_a\bar{\xi}_{\bar{b}})_i(x_0)=-c_{abij}\xi_a\bar{\xi}_{\bar{b}}(x_0)\\
&+\xi_{a,i\bar{j}}(x_0)\bar{\xi}_{\bar{a}}+\xi_a(x_0)\bar{\xi}_{\bar{a},i\bar{j}}(x_0)+\xi_{a,i}(x_0)\bar{\xi}_{\bar{a},\bar{j}}(x_0)+\xi_{a,\bar{j}}(x_0)\xi_{a,i}(x_0).
\end{split}
\end{equation}
In the above calculation,  we used that $g_{a\bar{b}}(x_0)=\delta_{ab}$,  $g_{a\bar{b},i}(x_0)=0$.  Next we need to estimate $D\xi$ and $D^2\xi$.
For this we have:
\begin{lem}
Let $\xi(z)$ be defined by $\exp_z(\xi(z))=\phi(z)$ with $\phi(z)=\exp_{x_0}(\xi_0)+N\cdot z$ in a neighborhood of $x_0$,  where we have taken normal coordinates centered at $x_0$.
Then there is a constant $C>0$ depending only on the background metric,  such that
\begin{equation*}
|D\xi|(x_0)\le C|\xi_0|_{x_0}^2,\,\,\,|D^2\xi|(x_0)\le C|\xi_0|_{x_0}^2.
\end{equation*}
\end{lem}
\begin{proof}
The result follows from differentiating the equality $\exp_z(\xi(z))=\phi(z)$.  Differentiate once,  we get:
\begin{equation}\label{2.17}
D_z\phi=D_z\exp_z(\xi)|_{\xi=\xi(z)}+D_{\xi}\exp_z(\xi)|_{\xi=\xi(z)}\cdot D_z\xi(z).
\end{equation}
Differentiate in $z$ again:
\begin{equation}\label{2.18}
\begin{split}
&D_z^2\phi=D_z^2\exp_z(\xi)|_{\xi=\xi(z)}+2D^2_{z\xi}\exp_z(\xi)|_{\xi=\xi(z)}\cdot D_z\xi(z)\\
&+(D_z\xi)^T\cdot D^2_{\xi\xi}\exp_z(\xi)|_{\xi=\xi(z)}\cdot D_z\xi+D_{\xi}\exp_z(\xi)\cdot D_z^2\xi(z).
\end{split}
\end{equation}
Now we evaluate (\ref{2.17}) at $z=0$,  so that we get:
\begin{equation*}
N=D_z\exp_z(\xi)|_{z=0,\xi=\xi_0}+D_{\xi}\exp_z(\xi)|_{z=0,\xi=\xi_0}D_z\xi(z)|_{z=0}.
\end{equation*}
From Lemma \ref{l2.19} below,  we have that:
\begin{equation*}
|D_z\exp_z(\xi)|_{z=0,\xi=\xi_0}-I|\le C|\xi_0|^2_{x_0},\,\,\,|D_{\xi}\exp_z(\xi)|_{z=0,\xi=\xi_0}-I|\le C|\xi_0|^2.
\end{equation*}
On the other hand,  we use (\ref{2.11N}),  and note that under normal coordinates, at $x_0$  we have:
\begin{equation*}
|g_{i\bar{j}}(w_0)-\delta_{ij}|\le Cd(w_0,x_0)^2\le C|\xi_0|^2_{x_0},
\end{equation*}
This would imply that we could choose $N$ so that $|N-I|\le C|\xi|_{x_0}^2$.
Hence
\begin{equation*}
D_z\xi|_{z=x_0}=(D_\xi\exp_z\xi)^{-1}|_{z=0,\xi=\xi_0}\big(N-D_z\exp_z(\xi)|_{z=0,\xi=\xi_0}\big).
\end{equation*}
Hence we can obtain that $|D_z\xi(x_0)|\le C|\xi_0|_{x_0}^2$.
To estimate $D_z^2\xi$,  we use (\ref{2.18}).  We note that,  $D_z^2\phi(z)=0$ and Lemma \ref{l2.19} gives $|D_z^2\exp_z(\xi)|_{z=0,\xi=\xi_0}|\le C|\xi_0|^2$.  Also we just need that $D_{z\xi}\exp_z(\xi)|_{z=0,\xi=\xi_0}$ and $D_{\xi\xi}^2\exp_z(\xi)|_{z=0,\xi=\xi_0}$ is bounded.
\end{proof}
In the above,  we used the following lemma about the Taylor expansion of the exponential map:
\begin{lem}\label{l2.19}
The following Taylor expansion holds for expential map on a K\"ahler manifold,  near $z=0$,  $\xi=0$ under normal coordinates centered at $z=0$:
\begin{equation*}
\exp_z(\xi)_m=z_m+\xi_m+\sum_{j,k,l}c_{jklm}(\frac{1}{2}\bar{z}_k+\frac{1}{6}\bar{\xi}_k)\xi_j\xi_l+O(|\xi|^2(|z|^2+|\xi|^2)).
\end{equation*}
\end{lem}
\begin{proof}
This is the equation (2.7) in \cite{DP} and a detailed proof can be found in Section 2 of \cite{DP}.
\end{proof}
Using (\ref{2.15N}),  we see that:
\begin{lem}\label{l2.20}
There exists a constant $C>0$ depending only on the background metric such that:
\begin{equation*}
|(|\xi(z)|^2_z)_{z_i\bar{z_j}}(x_0)+\sum_{a,b}c_{abij}\xi_{0,a}\bar{\xi}_{0,\bar{b}}|\le C|\xi_0|_{x_0}^3.
\end{equation*}
\end{lem}
Combining Lemma \ref{l2.10},  Proposition \ref{p2.14} and Lemma \ref{l2.20},  we conclude the following corollary:
\begin{cor}\label{c2.22N}
Let $\varphi\in C(M)$ and $\lambda[h_{\varphi}]\in\Gamma$ in the viscosity sense.  Define $\varphi_{\eps}$ according to (\ref{1.3N}).   Let $P\in C^2(M)$ and touch $\varphi_{\eps}$ at $x_0\in M$.  Assume that $M$ has nonnegative holomorphic bisection curvature at $x_0\in M$,  in the sense that $c_{ijkl}\eta_i\bar{\eta}_j\xi_k\bar{\xi}_l\ge 0$.  Then we have:
\begin{equation*}
\lambda\big(g^{i\bar{k}}(g_{j\bar{k}}+\partial_{z_j\bar{z}_k}P+C\eps^{\frac{1}{2}}g_{j\bar{k}})\big)(x_0)\in\Gamma.
\end{equation*}
Here the constant $C$ depends on the background metric and $||\varphi||_{L^{\infty}}$.

If in addition,  we assume that $\varphi\in C^{\alpha}$,  we can get:
\begin{equation*}
\lambda\big(g^{i\bar{k}}(g_{j\bar{k}}+\partial_{z_j\bar{z}_k}P+C\eps^{\frac{1+\alpha}{2-\alpha}}g_{j\bar{k}})\big)(x_0)\in\Gamma.
\end{equation*}
Here $C$ depends on the background metric and $C^{\alpha}$ norm of $\varphi$.
\end{cor}
\begin{proof}
First let us make the curvature assumption but do not assume $\varphi\in C^{\alpha}$.  Then Lemma \ref{l2.10} tells us that $|\xi_0|_{x_0}\le C\eps^{\frac{1}{2}}$.  Hence Proposition \ref{p2.14} would give us that:
\begin{equation*}
\lambda\big(g^{i\bar{k}}(g_{j\bar{k}}+C\eps^{\frac{1}{2}}g_{j\bar{k}}+\partial_{z_j\bar{z}_k}P-\sum_{a,b}c_{abjk}\xi_{0,a}\bar{\xi}_{0,\bar{b}})\big)(x_0)\in \Gamma.
\end{equation*}
The curvature assumption tells us that $\big(\sum_{a,b}c_{abij}\xi_{0,a}\bar{\xi}_{0,\bar{b}}\big)_{1\le i,j\le n}$ is semi-positive definite.  It follows that 
\begin{equation*}
\lambda\big(g^{i\bar{k}}(g_{j\bar{k}}+C\eps^{\frac{1}{2}}g_{j\bar{k}}+\partial_{z_j\bar{z}_k}P)\big)(x_0)\in\Gamma.
\end{equation*}

If we assume that $\varphi\in C^{\alpha}(M)$,  then we can estimate:
\begin{equation*}
\frac{1}{\eps}(|\xi(z)|^2_z)_{z_i\bar{z_j}}\le \frac{C|\xi_0|^3_{x_0}}{\eps}\delta_{ij}\le C\eps^{-1}(\eps[\varphi]_{\alpha})^{\frac{3}{2-\alpha}}\le C'\eps^{\frac{1+\alpha}{2-\alpha}}.
\end{equation*}
Here we used Lemma \ref{l2.10} and Lemma \ref{l2.20}.
\end{proof}
\begin{rem}
This is the only place in the whole proof that we use the curvature assumption.
\end{rem}

From now on,  let us estimate the $L^1$ difference between $\varphi_{\eps}$ and $\varphi$.
For this we need a version of mean value inequality on manifolds which is in light of \cite{C}.  
\begin{prop}\label{mean value inequality with a error term}
Let $u$ be a upper semicontinuous function defined on an open subset $U$ of a Riemannian manifold ($M,g$). Suppose that $\Delta_g u \ge 0$, $u$ is bounded and the Ricci curvature is bounded on $M$. Then there exists a constant $C$ which is a uniform constant depending on the bound of $||u||_{L^{\infty}}$ and $(M,g)$ such that:
\begin{equation*}
  u(x)\le \frac{1}{\alpha(2n)r^{2n}}\int_{B_r(x)}udvol_g+Cr,
\end{equation*}
for all $B_r(x)\subset U$.  Here $\alpha(2n)$ is the volume of the unit ball in $\mathbb{R}^{2n}$.
\end{prop}

\begin{proof}
Fix $r_0>0$ such that $B_{r_0}(x)\subset U$. Since $u$ is upper semicontinuous,  there exists a sequence of continuous functions $\{u_k\}$ such that $u\le u_k$ and $\lim_{k\rightarrow \infty}u_k=u$.  For example,  we could take $u_k(x)=\sup_{y\in M}\big(u(y)+\frac{1}{k}-kd_g^2(x,y)\big)$.  
Let us assume that $u<-1$ for the moment so that we can assume that each $u_k<0$.  

Let $0<r\le r_0$ and $h_k$ be the harmonic extension of $u_k$ on the $B_{r}(x)$.  By the maximum principle, we have that $0>h_k\ge u$.  This is because $h_k=u_k\ge u$ on $\partial B_{r}(x_0)$ and $\Delta_g(h_k-u)\le 0$ in $B_{r}(x_0)$.

By the Cheng-Yau harnack inequality (\cite{SY},  Chapter 1,  Section 3,  Theorem 3.1 ),  we have that:
\begin{equation*}
    |\nabla \log(-h_k)|_{B_{\frac{r}{2}}(x)}\le C_r,
\end{equation*}
where $C_r$ depends on the Ricci curvature and $r$.

Then we have that for any $s<\frac{r}{2}$,
\begin{equation*}
    \inf_{B_s(x)}h_k \ge e^{C_rs}\sup_{B_s(x)}h_k\ge e^{C_rs}\varphi(x).
\end{equation*}
Combining this with the almost monotonicity Lemma \ref{monotonicity lemma} we have that:
\begin{equation*}
\begin{split}
    r^{1-2n}\int_{\partial B_r(x)}h_k &\ge s^{1-2n}\int_{\partial B_s(x)}h_k -Cr\\
    &\ge e^{C_rs}s^{1-2n}vol(\partial B_s(x))\varphi(x)-Cr.
\end{split}
\end{equation*}
Hence, we have that:
\begin{equation*}
\begin{split}
    &\frac{1}{vol(\partial B_r(x))}\int_{\partial B_r(x)}u =\lim_{k\rightarrow \infty}\frac{1}{vol(\partial B_r(x))}\int_{\partial B_r(x)}h_k \\
    &\ge \lim_{s\rightarrow 0+}\frac{s^{1-2n}vol(\partial B_s(x))e^{C_rs}u(x)-Cr}{r^{1-2n}vol(\partial B_r(x))}\\
    &\ge \frac{2n\alpha(2n)r^{2n-1}}{vol(\partial B_r(x))}u(x)-Cr.
\end{split}
\end{equation*}
Here we use the Lemma \ref{estimate of the volume of balls} below.
So we have that,  for any $0<r\le r_0$:
\begin{equation*}
    2n \alpha(2n) r^{2n-1}u(x)\le \int_{\partial B_r(x)}u +Cr^{2n}.
\end{equation*}
Now we integrate in $r$,  we get:
\begin{equation}\label{2.25NNN}
\alpha(2n)r_0^{2n}u(x)\le \int_{B_{r_0}(x)}u-Cr_0^{2n+1}.
\end{equation}
So we have finished the proof if we assume that $\varphi\le -1$.  In general,  we consider $\tilde{u}:=u-(||u||_{L^{\infty}}+1)$ so that $\tilde{u}\le -1$.  Applying (\ref{2.25NNN}),  we see that
\begin{equation*}
\alpha(2n)r_0^{2n}\tilde{u}(x)\le \int_{B_{r_0}(x)}\tilde{u}dvol_g+Cr_0^{2n+1}.
\end{equation*}
Translating back to $u$,  we get:
\begin{equation*}
\begin{split}
&\alpha(2n)r_0^{2n}u(x)\le \int_{B_{r_0}(x)}udvol_g+(||u||_{L^{\infty}}+1)(\alpha(2n)r_0^{2n}-vol(B_{r_0}(x)))\\
&+Cr_0^{2n+1}\le C'r_0^{2n+1}.
\end{split}
\end{equation*}
\end{proof}

\begin{lem}\label{estimate of the volume of balls}
Let $(M,g)$ be a compact Riemannian manifold with real dimension $n$.       Then there exist $r_0>0$ and a constant $C$ depending only on the manifold such that for any $y\in M$ and $r<r_0$, we have the estimate:
        \begin{equation*}
\begin{split}
&|vol(\partial B_r(y))-n\alpha(n)r^{n-1}|\le Cr^n,\\
& |vol(B_r(y))-\alpha(n)r^{n}|\le C r^{n+1}.
\end{split}
        \end{equation*}
\end{lem}
\begin{proof}
We want to compare the volume form induced by the K\"ahler form with the Euclidean volume form in the normal coordinates.
Using the covering argument, it suffices to prove that for any $x\in M$, the conclusion of the lemma holds uniformly near $x$ with some $r_0$ and $C$. 

We just need to prove the first inequality,  as the second inequality would follow from the first by integrating in $r$.  First we find,  by using the Area formula:
\begin{equation*}
vol(\partial B_r(x))=\int_{S^{n-1}}|J_{\xi\mapsto \exp_x(r\xi)}|d\mathcal{H}^{n-1}(\xi).
\end{equation*}
In the above,  $\partial B_r$ is the image of the unit sphere (in $T_xM$) under the map $S^{n-1}\ni \xi\mapsto \exp_x\xi$.  $\mathcal{H}^{n-1}$ is the Hausdorff measure induced by the metric $g(x)$ on $T_xM$.  Also $J_{\xi\mapsto \exp_x(r\xi)}$ is the Jacobian of the map,  which is calculated by taking an othornormal frame of $S^{n-1}$ denoted as $v_1\wedge v_2\cdots\wedge v_{n-1}$,  then $J$ is the sup of $|(d_{\xi}(\exp(r\xi))v_1)\wedge \cdots \wedge (d_{\xi}(\exp(r\xi))v_{n-1}|$ over all orthonormal frames of $T_{\xi}S^{n-1}$.  Note that $d_{\xi}(\exp_z(\xi))|_{\xi=0}=id$,  which has Jacobian equaling $1$,  hence
\begin{equation*}
||d_{\zeta}\exp_z(\zeta)|_{\zeta=r\xi}-id||\le Cr,\,\,\xi\in S^{n-1}.
\end{equation*} hence we see that:
\begin{equation*}
J_{\xi\mapsto \exp_x(r\xi)}=r^{n-1}(1+O(r)).
\end{equation*}
Therefore,
\begin{equation*}
vol(\partial B_r(x))=\int_{S^{n-1}}r^{n-1}(1+O(1))d\mathcal{H}^{n-1}(\xi)=r^{n-1}vol(S^{n-1},g_E)+O(r^n).
\end{equation*}
Here $g_E$ just means the Euclidean metric.  But $vol(S^{n-1},g_E)=n\alpha(n)$.

\end{proof}

Another lemma that is used in the proof of Propositon \ref{mean value inequality with a error term} is:
\begin{lem}\label{monotonicity lemma}
Let $M$ be a closed manifold of dimension $2n$ with bounded Ricci curvature. Let $u$ be a bounded negative harmonic function defined on an open subset of $M$. Then we have that for $0<s<r$:
\begin{equation*}
    r^{1-2n}\int_{\partial B_r(x)}\varphi -s^{1-2n}\int_{\partial B_s(x)}\varphi \ge -Cr,
\end{equation*}
where $C$ is a uniform constant depending on $||\varphi||_{\infty}$ and $(M,\omega_0)$.
\end{lem}
\begin{proof}
Using the Laplacian comparision theorem, we have that $H\le H'$, where $H$ is the mean curvature of a geodesic sphere on $M$ and $H'$ is the mean curvature of a geodesic sphere on a space-form whose metric can be written in a normal coordinate as:
\begin{equation*}
    g_R=dr^2 + R^2 sinh^2(\frac{r}{R})d S_{2n-1}^2,
\end{equation*}
Here $R$ only depends on the lower bound of the Ricci curvature on $M$.
By calculation, we have that
\begin{equation*}
    H'(r)=\frac{(2n-1)cosh(\frac{r}{R})}{R sinh(\frac{r}{R})}=\frac{(2n-1)}{r}+O(1)
\end{equation*}
So we have that:
\begin{equation*}
\begin{split}
    \frac{d}{dt}(r^{1-2n}\int_{\partial B_r}\varphi)&=(1-2n)r^{-2n}\int_{\partial B_r}\varphi +r^{1-2n}\int_{\partial B_r}\varphi H \\
    &\le (1-2n)r^{-2n}\int_{\partial B_r}\varphi +r^{1-2n}\int_{\partial B_r}\varphi H' \\
    &\le (1-2n)r^{-2n}\int_{\partial B_r}\varphi +r^{1-2n}\int_{\partial B_r}\varphi(\frac{2n-1}{r}+O(1))\\
    &=O(1)
\end{split}
\end{equation*}
The conclusion of the lemma follows by integrating the above formula from $s$ to $r$.
\end{proof}

With the help of Proposition \ref{mean value inequality with a error term},  we are ready to estimate $\varphi_{\eps}-\varphi$ in $L^1$.
First,  we observe that by taking $\xi=0$:
\begin{equation*}
\varphi_{\eps}(z)=\sup_{\xi\in T_zM}\big(\varphi(\exp_z(\xi))+\eps-\frac{1}{\eps}|\xi|^2_z)\ge \varphi(z)+\eps.
\end{equation*}
Hence it will suffice to estimate $\int_M(\varphi_{\eps}-\varphi)\omega_0^n$.
We will do it in the following lemma:
\begin{lem}\label{l2.29}
Let $\varphi\in C(M)$ and satisfy $\lambda[h_{\varphi}]\in\Gamma$ in the viscosity sense.  Let $\varphi_{\eps}$ be defined as (\ref{1.3N}),  then one has:
\begin{equation*}
||\varphi_{\eps}-\varphi||_{L^1}\le C\eps^{\frac{1}{4}}.
\end{equation*}
Here $C$ depends only on the background metric and also $||\varphi||_{L^{\infty}}$.

If in addition,  we assume that $\varphi\in C^{\gamma}$,  then one has
\begin{equation*}
||\varphi_{\eps}-\varphi||_{L^1}\le C\eps^{\frac{1}{2(2-\gamma)}},
\end{equation*}
where $C$ depends only on the background metric and also $||\varphi||_{C^{\gamma}}.$
\end{lem}
\begin{proof}
We can find a cover $\{U_{\alpha}\}$ of $M$ such that $U_{\alpha}=B_{\sigma_0}(x_{\alpha})$ where $\sigma_0$ is selected such that $B_{2\sigma_0}(x_{\alpha})$ is contained in a normal coordinate neighborhood centered at $x_{\alpha}$.  Moreover,  there exists a K\"ahler potential $\phi_{\alpha}$ defined in $B_{2\sigma_0}(x_{\alpha})$. Namely $\omega_0=\sqrt{-1}\partial\bar{\partial}\rho_{\alpha}$ on $U_{\alpha}$.  
We can also assume that $\rho_{\alpha}\le 0$ and both $|\rho_{\alpha}|_{C^0}$ and $|\rho_{\alpha}|_{C^{\alpha}}$ are uniformly bounded. 

Since $\lambda[h_{\varphi}]\in\Gamma$ in the viscosity sense,  we see that $\lambda\big(g^{i\bar{k}}((\rho_{\alpha})_{j\bar{k}}+\varphi_{j\bar{k}})\big)\in\Gamma$ in the viscosity sense.  Since $\Gamma\subset \Gamma_+=\{\lambda\in\mathbb{R}^n:\sum_i\lambda_i>0\}$,  we see that
\begin{equation*}
\Delta(\varphi+\rho_{\alpha})>0\,\,\,\text{ on $U_{\alpha}$.}
\end{equation*}
Now we wish to apply the Proposition \ref{mean value inequality with a error term} to $\varphi+\rho_{\alpha}$.  For any $x\in M$,  we can find $\xi_x\in T_xM$ with $|\xi_x|_x\le (2||\varphi||_{L^{\infty}}\eps)^{\frac{1}{2}}$ by Lemma \ref{l2.10}.  Let $0<\beta<\frac{1}{2}$ be a constant we will choose later,  we have
\begin{equation}\label{uepsilon mean value inequality}
\begin{split}
    &\varphi_{\eps}(x) +\rho_{\alpha}(\exp_x(\xi_x))\le (\varphi+\rho_{\alpha})(\exp_x(\xi_x))+\eps \\
&\le \frac{1}{\alpha(2n)\eps^{2n\beta}}\int_{B_{\eps^{\beta}}(\exp_x(\xi_x))}(\varphi(y)+\rho_{\alpha}(y)) +C\eps^{\beta} \\
    &\le \frac{1}{\alpha(2n)\eps^{2n\beta}}\int_{B_{\eps^{\beta}-C_0\eps^{\frac{1}{2}}}(x)}\varphi(y) +\frac{1}{\alpha(2n)\eps^{2n\beta}}\int_{B_{\eps^{\beta}}(\exp_x(\xi_x))}\rho_{\alpha}(y)\\
    &-\rho_{\alpha}(\exp_x(\xi_x))+\rho_{\alpha}(\exp_x(\xi_x)) +C \eps^{\beta} \\
    &\le \frac{1}{\alpha(2n)\eps^{2n\beta}}\int_{B_{\eps^{\beta}-C_0\eps^{\frac{1}{2}}}(x)} \varphi(y)+\rho_{\alpha}(\exp_x(\xi_x))+C'\eps^{\beta}.
\end{split}
\end{equation}
In the above,  the constants $C$ and $C'$ depends only on $||\varphi||_{L^{\infty}}$ and the  background metric.

In the first line above,  we used the definition of $\varphi_{\eps}$ given by (\ref{1.3N}).
In the second line,  we used that $B_{\eps^{\beta}-C_0\eps^{\frac{1}{2}}}(x)\subset B_{\eps^{\beta}}(\exp_x(\xi_x))$,  which again follows from Lemma \ref{l2.10}.  Also $\varphi$ is negative.

In the last inequality,  we used that for any $y\in B_{\eps^{\beta}}(\exp_x(\xi_x))$,  $|\rho_{\alpha}(y)-\rho_{\alpha}(\exp_x(\xi_x))|\le \sup|\nabla\rho_{\alpha}|\eps^{\beta}$,  hence
\begin{equation*}
\begin{split}
&\frac{1}{\alpha(2n)\eps^{2n\beta}}\int_{B_{\eps^{\beta}}(\exp_x(\xi_x))}\rho_{\alpha}(y)-\rho_{\alpha}(\exp_x(\xi_x))\\
&=\frac{1}{\alpha(2n)\eps^{2n\beta}}\int_{B_{\eps^{\beta}}(\exp_x(\xi_x))}(\rho_{\alpha}(y)-\rho_{\alpha}(\exp_x(\xi_x)))+\big(\frac{vol(B_{\eps^{\beta}}(\exp_x(\xi_x)))}{\alpha(2n)\eps^{2n\beta}}-1\big)\rho_{\alpha}(\exp_x(\xi_x))\\
&\le C(|\rho_{\alpha}|_{C^1})\eps^{\beta}.
\end{split}
\end{equation*}
Here we also used Lemma \ref{estimate of the volume of balls}.

The result of (\ref{uepsilon mean value inequality}) is that:
\begin{equation}\label{2.31}
\varphi_{\eps}(x)\le \frac{1}{\alpha(2n)\eps^{2n\beta}}\int_{B_{\eps^{\beta}-C_0\eps^{\frac{1}{2}}}(x)}\varphi(y)\omega_0^n(y)+C\eps^{\beta}.
\end{equation}
Here $C$ depends only on the background metric and $||\varphi||_{L^{\infty}}$.
Now we can integrate on $M$ to get:
\begin{equation*}
\begin{split}
&\int_M\varphi_{\eps}(x)\omega_0^n(x)\le \frac{1}{\alpha(2n)\eps^{2n\beta}}\int_M\omega_0^n(x)\int_{B_{\eps^{\beta}-C_0\eps^{\frac{1}{2}}}(x)}\varphi(y)\omega_0^n(y)+Cvol(M)\eps^{\beta}\\
&=\int_M\varphi(y)\frac{vol(B_{\eps^{\beta}-C_0\eps^{\frac{1}{2}}}(y))}{\alpha(2n)\eps^{2n\beta}}\omega_0^n(y)+C_1\eps^{\beta}\le \inf_{y\in M}\frac{vol(B_{\eps^{\beta}-C_0\eps^{\frac{1}{2}}}(y))}{\alpha(2n)\eps^{2n\beta}}\int_M\varphi(y)\omega_0^n(y)+C_1\eps^{\beta}\\
&\le \frac{(\eps^{\beta}-C_0\eps^{\frac{1}{2}})^{2n}}{\eps^{2n\beta}}\int_M\varphi(y)\omega_0^n(y)+C_2\eps^{\beta}.
\end{split}
\end{equation*}
Therefore
\begin{equation*}
\int_M(\varphi_{\eps}-\varphi)\omega_0^n(x)\le \big(1-(1-C_0\eps^{\frac{1}{2}-\beta})^{2n}\big)\int_M\varphi(y)\omega_0^n(y)+C_2\eps^{\beta}\le C_3\eps^{\frac{1}{2}-\beta}+C_2\eps^{\beta}.
\end{equation*}
Hence we may take $\beta=\frac{1}{4}$ to obtain that $\int_M(\varphi_{\eps}-\varphi)\omega_0^n\le C\eps^{\frac{1}{4}}$.

Next if $\varphi\in C^{\gamma}$,  the calculation is very similar except that we have a better estimate for $\xi_x$: $|\xi_x|_x\le  C\eps^{\frac{1}{2-\gamma}}$.  
Then in \ref{uepsilon mean value inequality},  the range for $\beta$ should be $0<\beta<\frac{1}{2-\gamma}$.  Also instead of (\ref{2.31}),  one has instead:
\begin{equation*}
\varphi_{\eps}(x)\le \frac{1}{\alpha(2n)\eps^{2n\beta}}\int_{B_{\eps^{\beta}-C_1\eps^{\frac{1}{2-\gamma}}}(x)}\varphi(y)\omega_0^n(y)+C\eps^{\beta}.
\end{equation*}
Then we will make the choice that $\beta=\frac{1}{2(2-\gamma)}$ to obtain that:
\begin{equation*}
\int_M(\varphi_{\eps}-\varphi)\omega_0^n\le C\eps^{\frac{1}{2(2-\gamma)}}.
\end{equation*}
\end{proof}
Now we are in a position to prove Theorem \ref{uepsilon mean value inequality}.
\begin{proof}
(of Theorem \ref{uepsilon mean value inequality})
The proof follows from combining the lemmas obtained in this section.  Indeed,  point (i) follows from Lemma \ref{l2.11New}.  Point (ii) and (ii)' follow from Corollary \ref{c2.22N}.  Lemma \ref{l2.29} would give us point (iii) and (iii)'.
\end{proof}

\section{Stability estimates of Hessian equations}
\label{Stability}
Next we improve the estimate of $\varphi_{\eps}-\varphi$ from $L^1$ norm to the $L^{\infty}$ norm.
First we can prove the following proposition.
\begin{prop}\label{p3.2}
Let $v\in C(M)$,  $s>0$,  $\delta>0$and $\lambda\big(g^{i\bar{k}}((1+\frac{\delta}{2})g_{j\bar{k}}+v_{j\bar{k}})\big)_j^i\in\Gamma$ in the viscosity sense as in Definition \ref{d2.3}.  Define $A_{\delta,s}=\int_M ((1-\delta)v-\varphi-s)^+ e^{nF} \omega_0^n$.\\
Then there is a constant $\beta_n>0$, and $C>0$, depending only on the background metric and structural constant of f, such that:
\begin{equation*}
    \int_M exp(\frac{\beta_n(((1-\delta)v-\phi-s)^+)^{\frac{n+1}{n}}}{A_{\delta,s}^{\frac{1}{n}}})\omega_0^n \le C exp(C\delta^{-(n+1)}A_{\delta,s}).
\end{equation*}
\end{prop}
\begin{proof}
We consider the function 
\begin{equation*}
\Phi=\eps_1((1-\delta)v-\varphi-s)-(-\psi+\Lambda)^{\alpha}.
\end{equation*}
Here $\psi$ is defined as the solution to:
\begin{equation*}
\begin{split}
&(\omega_0+\sqrt{-1}\partial\bar{\partial}\psi)^n=\frac{((1-\delta)v-\varphi-s)^+e^{nF}}{A_{\delta,s}}\omega_0^n,\\
&\sup_M\psi=0.
\end{split}
\end{equation*}
Since $v$ is only continuous,  there is no guarantee that $\psi$ is $C^2$.  To do it rigorously,  we need to take $v_j$ to be a sequence of smooth functions which converges to $v$ uniformly,  and $\eta_j(x)$ be a sequence of smooth and strictly positive functions which converges to $x^+$.  Then we consider instead $\psi_j$,  which solves 
\begin{equation*}
(\omega_0+\sqrt{-1}\partial\bar{\partial}\psi_j)^n=\frac{\eta_j((1-\delta)v_j-\varphi-s)}{A_{j,\delta,s}}\omega_0^n.
\end{equation*}
Here $A_{j,\delta,s}=\int_M\eta_j((1-\delta)v_j-\varphi-s)^+e^{nF}\omega^n$.  Then we can consider $\Phi_j=\eps_0((1-\delta)v-\varphi-s)-(-\psi_j+\Lambda)^{\alpha}$.

For now,  let us assume that $\psi$ is $C^2$ so as to justify the following argument.  One can consider $\psi_j$ above in order to be completely rigorous (but argument is the same).

Since $\Phi$ is continuous,  it achieves maximum on $M$,  say,  at $p_0$,  and we denote:
\begin{equation*}
M_0:=\max_M\eps_0\big((1-\delta)v-\varphi-s)-(-\psi+\Lambda)^{\alpha_1}.
\end{equation*}
Our goal is to show that $M_0\le 0$.  If not,  namely $M_0>0$,  we can do the following argument.

At point $p$,  the function $Q:=\frac{1}{1-\delta}\big(\varphi+s+\frac{M_0+(-\psi+\Lambda)^{\frac{n}{n+1}}}{\eps_1}\big)$ touches $v$ from above at $p$.  Hence we get:
\begin{equation*}\label{0.18NNN}
\lambda((\frac{\delta}{2}+1)I+g^{i\bar{k}}Q_{j\bar{k}})\in\Gamma.
\end{equation*}
As before,  we consider the linearized operator defined as:
\begin{equation*}
Lv=\frac{\partial F}{\partial h_{ij}}((h_{\varphi})_j^i)g^{i\bar{k}}v_{j\bar{k}}.
\end{equation*}
We consider normal coordinate at $p$ so that $g_{i\bar{j}}=\delta_{ij}$ and $(h_{\varphi})_j^i$ is diagnal at $p$.  Then we have: $\frac{\partial F}{\partial h_{ij}}(p)=\frac{\partial f}{\partial \lambda_i}(\lambda(h_{\varphi}))\delta_{ij}$.  Then we can compute:
\begin{equation*}
L(\frac{Q}{1+\frac{\delta}{2}})=-\sum_i\frac{\partial f}{\partial \lambda_i}+\sum_i\frac{\partial f}{\partial \lambda_i}(\lambda(h_{\varphi}))(h_{\frac{Q}{1+\frac{\delta}{2}}})_{i\bar{i}}\ge -\sum_i\frac{\partial f}{\partial \lambda_i}+\sum_i\frac{\partial f}{\partial \lambda_i}(\lambda(h_{\varphi}))\mu_i.
\end{equation*}
In the above,  $\mu_i$ are the eigenvalues of $h_{\frac{Q}{1+\frac{\delta}{2}}}$,  in the increasing order.  The last inequality used Horn-Shur Lemma \ref{l3.2} (see \cite{Horn}).  Now we know that $\mu=(\mu_1,\cdots,\mu_n)\in\Gamma$,  thanks to (\ref{0.18NNN}).  
Hence 
\begin{equation*}
\sum_i\frac{\partial f}{\partial \lambda_i}(\lambda(h_{\varphi}))\mu_i=\frac{d}{dt}|_{t=0}(f(\lambda +t\mu)\ge 0.
\end{equation*}
It is due to that $t\mapsto f(\lambda+t\mu)$ is a concave function bounded from below.
Hence we get
\begin{equation*}
L(\frac{Q}{1+0.5\delta})\ge -\sum_i\frac{\partial f}{\partial \lambda_i}.
\end{equation*}
On the other hand,  we can compute:
\begin{equation*}
\begin{split}
&Q_{j\bar{k}}=\frac{1}{1-\delta}(\varphi_{j\bar{k}}+\eps_0^{-1}\frac{n}{n+1}(-\psi+\Lambda)^{-\frac{1}{n+1}}(-\psi_{j\bar{k}})-\frac{n}{(n+1)^2}(-\psi+\Lambda)^{-\frac{n+2}{n+1}}\psi_j\psi_{\bar{k}})\\
&=\frac{1}{1-\delta}\big((g_{\varphi})_{j\bar{k}}-\eps_0^{-1}\frac{n}{n+1}(-\psi+\Lambda)^{-\frac{1}{n+1}}(g_{\psi})_{j\bar{k}}-(1-\eps_0^{-1}\frac{n}{n+1}(-\psi+\Lambda)^{-\frac{1}{n+1}})g_{j\bar{k}}\\
&-\frac{n}{(n+1)^2}(-\psi+\Lambda)^{-\frac{n+2}{n+1}}\psi_j\psi_{\bar{k}})
\end{split}
\end{equation*}
So that (after dropping the term $\psi_j\psi_{\bar{k}}$ and replace $(-\psi+\Lambda)^{-\frac{1}{n+1}}$ by $\Lambda^{-\frac{1}{n+1}}$)
\begin{equation*}
\begin{split}
&-(1+0.5\delta)\sum_i\frac{\partial f}{\partial \lambda_i}\le LQ\le \frac{1}{1-\delta}\big(\sum_i\frac{\partial f}{\partial \lambda_i}\lambda_i-\eps_0^{-1}\frac{n}{n+1}(-\psi+\Lambda)^{-\frac{1}{n+1}}\frac{\partial F}{\partial h_{ij}}(g_{\psi})_{j\bar{i}}\\
&-\frac{1}{1-\delta}(1-\eps_0^{-1}\frac{n}{n+1}\Lambda^{-\frac{1}{n+1}})\sum_i\frac{\partial f}{\partial\lambda_i}.
\end{split}
\end{equation*}
Next,  you need to choose $\Lambda$ large enough so that 
\begin{equation*}
\eps_0^{-1}\frac{n}{n+1}\Lambda^{-\frac{1}{n+1}}=0.1\delta.
\end{equation*}
Then from above,  we get:
\begin{equation*}
\begin{split}
&0\le \sum_i\frac{\partial f}{\partial \lambda_i}\lambda_i-\eps_0^{-1}\frac{n}{n+1}(-\psi+\Lambda)^{-\frac{1}{n+1}}\frac{\partial F}{\partial h_{ij}}(g_{\psi})_{j\bar{i}}\\
&\le C_0f-\eps_0^{-1}\frac{n}{n+1}(-\psi+\Lambda)^{-\frac{1}{n+1}}\big(\frac{((1-\delta)v-u-s)^+e^{nF}}{A_{\delta,s}}\big)^{\frac{1}{n}}.
\end{split}
\end{equation*}
Since $M_0>0$,  it would imply that:
\begin{equation*}
(((1-\delta)v-\varphi-s)^+)^{\frac{1}{n}}(-\psi+\Lambda)^{-\frac{1}{n+1}}>\eps_0^{-\frac{1}{n}}.
\end{equation*}
Therefore,  if $\eps_0$ is small enough so that
\begin{equation*}
C_0-\eps_0^{-\frac{n+1}{n}}A_{\delta,s}^{-\frac{1}{n}}<0,
\end{equation*}
it would not happen.
\end{proof}
In the above,  we used the following Horn-Shur Lemma:
\begin{lem}\label{l3.2}
(\cite{Horn})
Let $A$ be a Hermitian matrix,  then the vector $(a_{1\bar{1}},a_{2\bar{2}},\cdots,a_{n\bar{n}})$ formed by the diagnal entries is contained in the convex envelope of $\{(\mu_{\sigma(1)},\mu_{\sigma(2)},\cdots,\mu_{\sigma(n)})\}_{\sigma\in S_n}$.  Here $\mu_i$ are the eigenvalues of $A$ and $\sigma$ denotes permutations.
\end{lem}

Next we need to estimate the upper level set $\{(1-\delta)v-\varphi-s>0\}$,  to the effect that De Giorge's lemma applies.  For this we have the following lemma:
\begin{lem}\label{l3.3}
Let $0<\delta<1$ and $v$ be as assumed in Proposition \ref{p3.2} and $\varphi$ solves $f(\lambda(h_{\varphi}))=e^F$ with $e^{nF}\in L^{p_0}(\omega_0^n)$ for some $p_0>1$.  Define $\Omega_{\delta,s}=\{(1-\delta)v-\varphi-s>0\}$.  Assume that $\delta$ and $s$ are chosen so that:
\begin{equation*}
\int_M((1-\delta)v-\varphi-s)^+e^{nF}\omega_0^n\le \delta^{n+1}.
\end{equation*}
Then for any $r>0$ and any $\mu<\frac{1}{nq_0}$,  we have:
\begin{equation*}
rvol(\Omega_{\delta,s+r})\le B_0vol(\Omega_{\delta,s})^{1+\mu}.
\end{equation*}
Here $q_0$ is the dual exponent to $p_0$,  namely $\frac{1}{q_0}=1-\frac{1}{p_0}$.  Here $B_0$ depends only on the background metric,  $p_0$ and choice of $\mu$.
\end{lem}
\begin{proof}
we see  Proposition \ref{p3.2} implies that for any $p\ge 1$,  we have:
\begin{equation*}
\int_{\Omega_{\delta,s}}((1-\delta)v-\varphi-s)^{\frac{n+1}{n}p}\omega_0^n\le C_1e^{C_0\delta^{-(n+1)}A_{\delta,s}}A_{\delta,s}^{\frac{p}{n}}.
\end{equation*}
In the above,  $C_1=C_0p!$.
If we put $p'=\frac{n+1}{n}p$,  we get that
\begin{equation}\label{0.0.5}
||((1-\delta)v-\varphi-s)^+||_{L^{p}}\le C_2e^{C_0\delta^{-(n+1)}A_{\delta,s}}A_{\delta,s}^{\frac{1}{n+1}}.
\end{equation}
Here $C_2$ depends on $p$ but $C_0$ doesn't.

Denote $h=((1-\delta)v-\varphi-s)^+$ for the moment. 
Now we use the assumption that $e^{nF}\in L^{p_0}(\omega_0^n)$ for some $p_0>1$ to obtain that:
\begin{equation*}
\begin{split}
&A_{\delta,s}=\int_{\Omega_{\delta,s}}he^{nF}\omega_0^n\le \big(\int_{\Omega_s}h^{q_0}\omega_0^n\big)^{\frac{1}{q_0}}||e^{nF}||_{L^{p_0}(\omega_0^n)}\le ||h||_{L^{\beta q_0}(\omega_0^n)}vol(\Omega_{\delta,s})^{\frac{1}{q_0}(1-\frac{1}{\beta})}||e^{nF}||_{L^{p_0}(\omega_0^n)}\\
&\le C_2e^{C_0\delta^{-(n+1)}A_{\delta,s}}A_{\delta,s}^{\frac{1}{n+1}}vol(\Omega_{\delta,s})^{\frac{1}{q_0}(1-\frac{1}{\beta})}||e^{nF}||_{L^{p_0}(\omega_0^n)}.
\end{split}
\end{equation*}
That is,  we have
\begin{equation*}
A_{\delta,s}^{\frac{1}{n+1}}\le C_2^{\frac{1}{n}}e^{\frac{C_0}{n}\delta^{-(n+1)}A_{\delta,s}}vol(\Omega_{\delta,s})^{\frac{1}{nq_0}(1-\frac{1}{\beta})}||e^{nF}||_{L^{p_0}(\omega_0^n)}^{\frac{1}{n}}.
\end{equation*}
On the other hand,  
\begin{equation*}
\begin{split}
&||h||_{L^1(\omega_0^n)}\le ||h||_{L^{\beta}(\omega_0^n)}vol(\Omega_{\delta,s})^{1-\frac{1}{\beta}}\le C_2e^{C_0\delta^{-(n+1)}A_{\delta,s}}A_{\delta,s}^{\frac{1}{n+1}}vol(\Omega_{\delta,s})^{1-\frac{1}{\beta}}\\
&\le C_2^{1+\frac{1}{n}}e^{(1+\frac{1}{n})C_0\delta^{-(n+1)}A_{\delta,s}}vol(\Omega_{\delta,s})^{(1+\frac{1}{nq_0})(1-\frac{1}{\beta})}||e^{nF}||_{L^{p_0}(\omega_0^n)}.
\end{split}
\end{equation*}
Note that 
\begin{equation*}
||h||_{L^1(\omega_0^n)}\ge rvol(\Omega_{\delta,s+r}).
\end{equation*}
So that we have
\begin{equation*}
rvol(\Omega_{\delta,s+r})\le C_3vol(\Omega_{\delta,s})^{(1+\frac{1}{nq_0})(1-\frac{1}{\beta})}.
\end{equation*}
Here for any $\mu<\frac{1}{nq_0}$,  one could take $\beta$ so that $(1+\frac{1}{nq_0})(1-\frac{1}{\beta})=1+\mu$.
\end{proof}
On the other hand,  we have the following lemma of De Giorge:
\begin{lem}\label{l3.4}
Let $\phi:[0,\infty)\rightarrow [0,\infty)$ be an decreasing and continuous function,  such that there exists $\mu>0$,  $B_0>0$,  $s_0\ge 0$,  such that for any $r>0$,  and any $s\ge s_0$,  one has:
\begin{equation*}
r\phi(s+r)\le B_0\phi(s)^{1+\mu}.
\end{equation*}
Then $\phi(s)\equiv 0$ for $s\ge s_0+\frac{2B_0\phi(s_0)^{\mu}}{1-2^{-\mu}}.$
\end{lem}
\begin{proof}
We can choose a sequence $\{s_k\}_{k\ge 1}$ by induction:
\begin{equation*}
s_{k+1}-s_k=2B_0\phi(s_k)^{\mu}.
\end{equation*}
Then we choose $s=s_k$,  $r=s_{k+1}-s_k$,  we see that:
\begin{equation*}
\phi(s_{k+1})\le \frac{B_0\phi(s_k)^{1+\mu}}{s_{k+1}-s_k}\le \frac{1}{2}\phi(s_k).
\end{equation*}
That is,  $\phi(s_k)\le 2^{-k}\phi(s_0)$.  Hence
\begin{equation*}
0\le s_{k+1}-s_k\le 2B_0\phi(s_0)^{\mu}2^{-k\mu}.
\end{equation*}
This implies that $s_k$ is increasing and bounded from above,  hence must converge to some $s_{\infty}$.  Moreover,
\begin{equation*}
s_{\infty}-s_0=\sum_{k=0}^{\infty}(s_{k+1}-s_k)\le \sum_{k=0}^{\infty}2B_0\phi(s_0)^{\mu}2^{-k\mu}=\frac{2B_0\phi(s_0)^{\mu}}{1-2^{-\mu}}.
\end{equation*}
\end{proof}
Combining Lemma \ref{l3.3} and Lemma \ref{l3.4},  we immediately have the following corollary:
\begin{cor}\label{c3.5}
Let $v$ and $\varphi$ be as in previous lemma.  Let $\delta>0$,  $s_0>0$ be chosen so that:
\begin{equation*}
A_{\delta,s_0}=\int_M((1-\delta)v-\varphi-s_0)^+e^{nF}\omega_0^n\le \delta^{n+1},
\end{equation*}
then for any $\mu<\frac{1}{nq_0}$, 
\begin{equation*}
\sup_M((1-\delta)v-\varphi)\le s_0+C_0vol(\{1-\delta)v-\varphi>s_0\})^{\mu}.
\end{equation*}
Here $C_0$ depends on $\mu$,  $||e^{nF}||_{L^{p_0}}$ but not on $\delta$.
\end{cor}
\begin{proof}
Put $\phi(s)=vol(\Omega_{\delta,s})$.  Then $\phi(s)$ would satisfy the assumptions of Lemma \ref{l3.4},  according to Lemma \ref{l3.3}.  Then we see that $\phi(s)\equiv 0$ for $s\ge s_0+\frac{2B_0\phi(s_0)^{\mu}}{1-2^{-\mu}}$.  This means $(1-\delta)v-\varphi\le s_0+\frac{2B_0\phi(s_0)^{\mu}}{1-2^{-\mu}}$.
\end{proof}
Next we need to get rid of $\delta$ in the above estimate.  For this let us assume that $v\le 0$ and is bounded.  Then in the left hand side,  we have
\begin{equation*}
\sup_M(v-\varphi)\le \sup_M((1-\delta)v-\varphi).
\end{equation*}
On the other hand,  we note that:
\begin{equation*}
vol(\Omega_{\delta,s_0})\le \frac{1}{s_0}\int_{\Omega_{\delta,s_0}}((1-\delta)v-\varphi)^+\omega_0^n\le \frac{1}{s_0}\big(||(v-\varphi)^+||_{L^1}+\delta||v||_{L^{\infty}}vol(\Omega_{\delta,s_0}).
\end{equation*}
Therefore,  if $\delta||v||_{L^{\infty}}\le \frac{1}{2}$,  we immediately conclude that:
\begin{equation}\label{3.6N}
vol(\Omega_{\delta,s_0})\le \frac{2}{s_0}||(v-\varphi)^+||_{L^1}.
\end{equation}
In view of (\ref{3.6N}) and Corollary \ref{c3.5},  we conclude that:
\begin{lem}\label{l3.7}
Let $v$ and $\varphi$ be as in Proposition \ref{p3.2}.  We assume additionally that $v\le 0$ and is bounded.  Assume also that $s_0$ and $\delta$ are chosen so that:
\begin{enumerate}
\item $s_0\ge 2\delta||v||_{L^{\infty}}$,
\item $A_{\delta,s_0}\le \delta^{n+1}$,  
\end{enumerate}
then for any $\mu<\frac{1}{q_0n}$,  we have:
\begin{equation}\label{3.8N}
\sup_M(v-\varphi)\le s_0+C_1s_0^{-\mu}||(v-\varphi)^+||_{L^1}^{\mu}.
\end{equation}
\end{lem}
Define $s_*(\delta)$ to be the infimum (with fixed $\delta$) of $s_0$ which satisfies the assumptions (1) and (2) put forth in Lemma \ref{l3.7}.
Note that the left hand side of (\ref{3.8N}) does not depend on $s_0$,  we get:
\begin{equation*}
\sup_M(v-\varphi)\le \inf_{s_0\ge s_*(\delta)}(s_0+C_1s_0^{-\mu}||(v-\varphi)^+||_{L^1}^{\mu}).
\end{equation*}
At this point,  it is clear that we must have an estimate (from above) of $s_*(\delta)$,  and we have:
\begin{lem}\label{l3.9}
Let $\delta>0$ and $v$ be as in Proposition \ref{p3.2}.  We also assume that $v$ is bounded.  Define $s_*(\delta)$ to be the infimum of the set of $s_0$ such that $A_{\delta,s_0}\le \delta^{n+1}$ and $s_0\ge 2\delta||v||_{L^{\infty}}$,  then for any $\beta>1$,
\begin{equation*}
s_*(\delta)\le \max(2\delta||v||_{L^{\infty}},C_2\delta^{-\frac{q_0n}{1-\beta^{-1}}}||(v-\varphi)^+||_{L^1}\big).
\end{equation*}
Here $C_2$ depends on $||e^{nF}||_{L^{p_0}}$,  $\beta$,  but is independent of $\delta$.
\end{lem}
\begin{proof}
Note that $A_{\delta,s_0}$ depends continuously on $s_0$,  we see that,  for $s_0=s_*(\delta)$,  either you have $A_{\delta,s_0}=\delta^{n+1}$ or $s_0=2\delta||v||_{L^{\infty}}$.

If $s_*(\delta)=2\delta||v||_{L^{\infty}}$,  then we are already done.

If $A_{\delta,s_*(\delta)}=\delta^{n+1}$,  then from Proposition \ref{p3.2},  we see that for any $q>1$,
\begin{equation*}
A_{\delta,s_*}^{-\frac{q}{n}}\int_M(((1-\delta)v-\varphi-s)^+)^{\frac{n+1}{n}q}\omega_0^n\le C\int_M\exp\big(\frac{\beta_n((1-\delta)v-\varphi-s_*)^{\frac{n+1}{n}}}{A_{\delta,s_*}^{\frac{1}{n}}}\big)\omega_0^n\le C'.
\end{equation*}
This gives:
\begin{equation}\label{3.10}
||((1-\delta)v-\varphi-s_*)^+||_{L^q}\le C_4 A_{\delta,s_*}^{\frac{1}{n+1}}.
\end{equation}
Here $C_4$ depends on $q$ but is independent of $\delta$.
Therefore
\begin{equation*}
\begin{split}
&A_{\delta,s_*}=\int_{\Omega_{\delta,s_*}}((1-\delta)v-\varphi-s)^+e^{nF}\omega_0^n\le ||e^{nF}||_{L^{p_0}}||((1-\delta)v-\varphi-s_*)^+||_{L^{q_0}}\\
&\le ||e^{nF}||_{L^{p_0}}||((1-\delta)v-\varphi-s_*)^+||_{L^{\beta q_0}}vol(\Omega_{\delta,s_*})^{\frac{1}{q_0}(1-\beta^{-1})}\\
&\le ||e^{nF}||_{L^{p_0}}C(\beta,p_0)A_{\delta,s_*}^{\frac{1}{n+1}}(\frac{2}{s_*}||(v-\varphi)^+||_{L^1})^{\frac{1}{q_0}(1-\beta^{-1})}.
\end{split}
\end{equation*}
In the above,  we used (\ref{3.6N}) and (\ref{3.10}).  Keeping in mind that $A_{\delta,s_*}=\delta^{n+1}$,  we get that
\begin{equation*}
s_*\le C\delta^{-\frac{nq_0}{1-\beta^{-1}}}||(v-\varphi)^+||_{L^1}.
\end{equation*}
\end{proof}
As a corollary of the above estimate,  we see that:
\begin{cor}\label{stability}
Let $\varphi$ be a bounded negative solution to $f(\lambda(h_{\varphi}))=e^{F}$ with $e^{nF}\in L^{p_0}(\omega_0^n)$ for some $p_0>1$.
Let $\{\varphi_{\eps}\}_{0<\eps<1}$ be a family of upper approximations of $\varphi$ such that $\varphi_{\eps}<0$,  $\varphi_{\eps}$ uniformly bounded,  and for some $\gamma_2>0$,  $0<\gamma_1<\frac{\gamma_2}{q_0n}$,
\begin{enumerate}
\item There exists $C_{3}>0$ independent of $\eps$ such that $\frac{C_{1.5}\eps^{\gamma_1}}{2}I+g^{i\bar{k}}(g_{j\bar{k}}+(\varphi_{\eps})_{j\bar{k}})\in \Gamma$ in the viscosity sense;
\item There exists $C_4>0$ independent of $\eps$ such that $||\varphi_{\eps}-\varphi||_{L^1}\le C\eps^{\gamma_2}.$
\end{enumerate} 
Then for each $\gamma<\min(\gamma_1,\gamma_2-q_0n\gamma_1,\frac{\gamma_2}{1+q_0n})$,  there exists $C_5>0$ such that
\begin{equation*}
\sup_M(\varphi_{\eps}-\varphi)\le C_5\eps^{\gamma}.
\end{equation*}
\end{cor}
\begin{proof}
Apply Lemma \ref{l3.7} with $\delta=C_3\eps^{\gamma_1}$,  we find that
\begin{equation}\label{3.12NN}
\sup_M(\varphi_{\eps}-\varphi)\le s_0+C_1s_0^{-\mu}(C_4\eps^{\gamma_2})^{\mu},
\end{equation}
for any $s\ge s_*(\delta)$.  (Of course here we have $\delta=C_{3}\eps^{\gamma_1}$.)
We would like to obtain an optimal estimate for $\sup_M(\varphi_{\eps}-\varphi)$ by minimizing the right hand side of (\ref{3.12NN}) in terms of $s_0$,  with $s_0$ subject to the constraint that $s_0\ge s_*(C_3\eps^{\gamma_1})$.

Denote $M_1=\sup_{\eps>0}||\varphi_{\eps}||_{L^{\infty}}$.  
Lemma \ref{l3.9} gives us that:
\begin{equation}
\begin{split}
&s_*(C_3\eps^{\gamma_1})\le \max(2(C_3\eps^{\gamma_1})M_1,
C_2(C_3\eps^{\gamma_1})^{-\frac{q_0n}{1-\beta^{-1}}}C_4\eps^{\gamma_2})\\
&\le C_{4.5}\eps^{\min(\gamma_1,-\frac{q_0n\gamma_1}{1-\beta^{-1}}+\gamma_2)},\,\,\,C_{4.5}=\max(2C_3M_1,\,C_2C_3^{-\frac{q_0n}{1-\beta^{-1}}}C_4).
\end{split}
\end{equation}
On the other hand,  the right hand side of (\ref{3.12NN}) is monotone decreasing for $s_0\in (0,(C_1C_4^{\mu})^{\frac{1}{1+\mu}}\eps^{\frac{\gamma_2\mu}{1+\mu}})$ and is monotone increasing for \sloppy $s_0\in ((C_1C_4^{\mu})^{\frac{1}{1+\mu}}\eps^{\frac{\gamma_2\mu}{1+\mu}},+\infty)$.
Denote $\gamma_3=\min(\gamma_1,-\frac{q_0n\gamma_1}{1-\beta^{-1}}+\gamma_2)$ and $C_{4.7}=(C_1C_4^{\mu})^{\frac{1}{1+\mu}}$.  There are two cases to consider:
\begin{enumerate}
\item If $\gamma_3>\frac{\gamma_2\mu}{1+\mu}$,  then the choice of $s_0=C_{4.7}\eps^{\frac{\gamma_2\mu}{1+\mu}}$ will be $\ge s_*(C_3\eps^{\gamma_1})$,  and we are allowed to use (\ref{3.12NN}) to conclude that
\begin{equation*}
\sup_M(\varphi_{\eps}-\varphi)\le C_{4.9}\eps^{\frac{\gamma_2\mu}{1+\mu}}.
\end{equation*}
\item If $\gamma_3\le \frac{\gamma_2\mu}{1+\mu}$,  then we may not be able to choose $s_0=C_{4.7}\eps^{\frac{\gamma_2\mu}{1+\mu}}$.  In this case,  we just choose $s_0=C_{4.5}\eps^{\gamma_3}$ in the estimate (\ref{3.12NN}) to conclude that:
\begin{equation*}
\sup_M(\varphi_{\eps}-\varphi)\le C_{4.91}\eps^{\gamma_3}.
\end{equation*}
\end{enumerate}
Combining both cases,  we get
\begin{equation*}
\sup_M(\varphi_{\eps}-\varphi)\le C_{4.92}\eps^{\min(\gamma_3,\frac{\gamma_2\mu}{1+\mu})}=C_{4.92}\eps^{\min(\gamma_1,\gamma_2-\frac{q_0n\gamma_1}{1-\beta^{-1}},\frac{\gamma_2\mu}{1+\mu})}.
\end{equation*}
Note that one can choose arbitrary $\beta>1$ and arbitrary $\mu<\frac{1}{q_0n}$,   so the H\"older exponent can be made as close to $\min(\gamma_1,\gamma_2-q_0n\gamma_1,\frac{\gamma_2}{1+q_0n})$ as one desires.
\end{proof}

We can now prove the H\"older continuity of solutions:
\begin{thm}\label{t4.1}
Let $(M,\omega_0)$ be a compact K\"ahler manifold with nonnegative holomorphic bisectional curvature.   Let $\varphi$ be a solution to the Hessian equation:
\begin{equation*}
f(\lambda[h_{\varphi}])=e^F,\,\,\,\lambda[h_{\varphi}]\in \Gamma,\,\,\,\sup_M\varphi=-1,
\end{equation*}
where $f$ and $\Gamma$ satisfies the structural assumptions put forth after (\ref{1.1New}).  Assume that $e^{nF}\in L^{p_0}(\omega_0^n)$ for some $p_0>1$.  Then for any $\mu<\frac{1}{2(1+q_0n)}$,  $||\varphi||_{C^{\mu}}\le C$.  Here $q_0=\frac{p_0}{p_0-1}$ and $C$ depends on the structural constants of $f$,  $||e^F||_{L^{p_0n}}$ and the background metric.
\end{thm}
\begin{proof}
Let $\frac{1}{2}<\beta\le 1$,  then we have:
\begin{equation*}
\varphi_{\eps}(z)\ge\sup_{|\xi|_z<\eps^{\beta}}\big(\varphi(\exp_z(\xi))+\eps-\frac{1}{\eps}|\xi|_z^2\big)\ge \sup_{d(w,z)<\eps^{\beta}}\varphi(w)+\eps-\eps^{2\beta-1}.
\end{equation*}
On the other hand,  using \ref{stability},  we see that 
\begin{equation*}
\sup_{d(w,z)<\eps^{\beta}}\varphi(w)-\varphi(z)\le \eps^{2\beta-1}+C\eps^{\gamma}.
\end{equation*}
Hence we may take $\beta=\frac{1+\gamma}{2}$,  we see that:
\begin{equation*}
\sup_{d(w,z)<\eps^{\frac{1+\gamma}{2}}}\varphi(w)-\varphi(z)\le C\eps^{\gamma}=C(\eps^{\frac{1+\gamma}{2}})^{\frac{2\gamma}{1+\gamma}}.
\end{equation*}
This will imply that $\varphi$ is H\"older continuous with exponent $\frac{2\gamma}{1+\gamma}$.  
According to Corollary \ref{stability},  we can take any $\gamma<\min(\gamma_1,\gamma_2-q_0n\gamma_1,\frac{\gamma_2}{1+q_0n})$.  
On the other hand,  from Theorem \ref{hessian L1 estimate},  we see that we could take $\gamma_1=\frac{1}{2}$,  $\gamma_2=\frac{1}{4}$ (with no initial assumption on the continuity of $\varphi$).  So the range of $\gamma$ given by Corollary \ref{stability} says:
\begin{equation*}
\gamma<\sup_{0<\gamma_1\le \frac{1}{2}}\min(\gamma_1,\frac{1}{4}-q_0n\gamma_1,\frac{1}{4(1+q_0n)}).
\end{equation*}
We see that the range of $\gamma$ is $\gamma<\frac{1}{4(1+q_0n)}$.  
Therefore the range of H\"older exponent is $<\frac{2}{4(1+q_0n)+1}$.

Next we observe that this H\"older exponent can be improved by an iteration process,  using the second part of Theorem \ref{hessian L1 estimate}.  To be more specific,  we already obtained some H\"older exponent $\mu_0$,  then the second part of Theorem \ref{hessian L1 estimate} would give an improved estimate for both $\gamma_1$ and $\gamma_2$ in Corollary \ref{stability},  hence will give us an improved $\gamma$,  hence an improved H\"older exponent.  Now we do the detailed calculation.

We already found that $\varphi$ has H\"older continuity for all $\mu<\frac{2}{4(1+q_0n)+1}$.  We denote $\mu_0=\frac{2}{4(1+q_0n)+1}$.  In general,  let us denote $0<\mu_k<1$ to be the upper bound of the H\"older exponent.   Then (ii)' of Theorem \ref{hessian L1 estimate} gives us that $\gamma_1$ can be improved to be $\frac{1+\mu_k}{2-\mu_k}$ and part (iii)' of the same theorem implies that $\gamma_2$ can be improved to be $\frac{1}{2(2-\mu_k)}$.  So the upper bound of $\gamma$ given by Corollary \ref{stability} is given by:
\begin{equation*}
\gamma<\sup_{0<\gamma_1\le \frac{1+\mu_k}{2-\mu_k}}\min(\gamma_1,\frac{1}{2(2-\mu_k)}-q_0n\gamma_1,\frac{1}{2(2-\mu_k)(1+q_0n)}).
\end{equation*}
Hence the upper bound of $\gamma$ is just $\frac{1}{2(2-\mu_k)(1+q_0n)}$.  Therefore,  the new upper bound of H\"older exponent can be calculated by:
\begin{equation*}
\mu_{k+1}=\frac{2\times \frac{1}{2(2-\mu_k)(1+q_0n)}}{1+\frac{1}{2(2-\mu_k)(1+q_0n)}}=\frac{2}{2(2-\mu_k)(1+q_0n)+1},\,\,\,\mu_0=\frac{2}{4(1+q_0n)+1}.
\end{equation*}
One can verify that $\mu_k\le \mu_{k+1}$ and $\mu_k$ converges to $\frac{1}{2(1+q_0n)}$.
\end{proof}

\end{document}